\documentclass[11pt]{article}
\usepackage{amsmath,amsthm,amsfonts,latexsym,amscd,amssymb}

\scrollmode

\swapnumbers
\theoremstyle{plain}
\newtheorem{theorem}{Theorem}[section]
\newtheorem{lemma}[theorem]{Lemma}
\newtheorem{corollary}[theorem]{Corollary}
\newtheorem{proposition}[theorem]{Proposition}

\theoremstyle{definition}
\newtheorem{definition}[theorem]{Definition}

\newtheorem{convention}[theorem]{Convention}

\theoremstyle{remark}
\newtheorem{remark}[theorem]{Remark}

\newcommand{\R}{\mathbb{R}}
\newcommand{\reals}{\mathbb{R}}

\newcommand{\N}{\mathbb{N}}
\newcommand{\Z}{\mathbb{Z}}
\newcommand{\integers}{\mathbb{Z}}

\newcommand{\obundle}{\widetilde{o}}  % orientation bundle for manifolds
\DeclareMathOperator{\id}{id}

\newcommand{\tensor}{\otimes}

\newcommand{\iso}{\cong}

\DeclareMathOperator{\Ker}{Ker}    %Kernels
\DeclareMathOperator{\Image}{Im}    %Images

\DeclareMathOperator{\scal}{scal}  % scalar curvature

\newcommand{\Kommentar}[1]{\textbf{(Kommentar: #1)}}
\newcommand{\forget}[1]{}
\newcommand{\Cliff}{{\mathcal C}\ell}
{\catcode`@=11\global\let\c@equation=\c@theorem}

\renewcommand{\theenumi}{\@arabic{enumi}}

\begin{document}
\date{Last edited: July 8, 99 --- last compiled: \today}

\title{Positive and negative results concerning the Gromov-Lawson-Rosenberg
  conjecture}
\author{ Michael
Joachim\thanks{\noindent email:
michael.joachim@math.uni-muenster.de\protect\\
www: ~http://www.math.uni-muenster.de/u/lueck/org/staff/joachim/\protect}\\ Fachbereich Mathematik\\ Universit\"at M\"unster\\
Einsteinstr.~62\\         48149 M\"unster \and
Thomas Schick\thanks{\noindent e-mail:
thomas.schick@math.psu.edu\protect\\
www:~http://www.math.psu.edu/schick/\protect\\ Stay at Penn state funded by the DAAD}\\ Dept.~of
Mathematics\\ Penn State University\\ 218 McAllister Building\\
University Park, PA 16802}
\maketitle

\begin{abstract}
The Gromov-Lawson-Rosenberg (GLR)-conjecture for a group $\Gamma$ 
states that a closed spin manifold $M^n$ ($n\ge5$) 
with fundamental group $ \Gamma$ admits a 
metric with $\scal>0$ if and only if its $C^*$-index $\alpha(M)\in
KO_n(C^*_{red}(\Gamma))$ vanishes.
We prove this for groups $\Gamma$ with low-dimensional classifying space and
products of such groups with free abelian groups, provided the
assembly map for the group $\Gamma$ is (split) injective (and $n$ large enough).

\forget{
We show that if for an $N$-dimensional group $\Gamma$ 
 the injectivity part of the Baum-Connes conjecture is
true (i.e.~the assembly map $A:KO_*(B\Gamma)\to KO_*(C^*_{red}\Gamma)$
is injective) then the unstable Gromov-Lawson-Rosenberg (GLR)-conjecture 
holds in dimensions $n\ge \max\{N-4,5\}$. That means  
a closed spin manifold $M^n$  with $\pi_1(M)\cong\Gamma$  admits a metric with
$\scal>0$ if and only if its $C^*$-valued index $\alpha(M)\in
KO_n(C^*_{red}\Gamma)$ vanishes. 
E.g.~it is true in full generality for hyperbolic groups of dimension at most $9$.
 In addition we show that if $\dim B\Gamma\le 4$ and
the assembly map is split injective for $\Gamma$ then the
GLR-conjecture holds for fundamental groups
isomorphic to $\Gamma\times\Z^k$ for all $k\in\N$.
}%end{forget}

On the other hand, we construct a
$5$-dimensional  spin
manifold $M$ which does not admit a metric with $\scal>0$ but
has the property 
that already the image of its
$KO$-orientation $pD[M]\in KO_*(B\pi_1(M))$ vanishes.
Therefore a corresponding weakened version
of the GLR-conjecture is wrong. 

Last we address non-orientable manifolds. We give a
reformulation of the minimal surface method of Schoen and Yau
(extended to dimension $8$) and introduce a
non-orientable (twisted) version of it. We then construct
a $5$-dimensional manifold whose orientation cover admits 
a metric of positive scalar curvature and use the latter to 
show that the manifold itself does not. The manifold also is
 a counterexample to a twisted analog of the
GLR-conjecture because its twisted index vanishes.

MSC-number: 53C20 (global Riemannian Geometry)
\end{abstract}

\section{Introduction}
This paper studies the question which closed manifolds admit a
Riemannian metric with positive scalar curvature. (Throughout the
paper, all manifolds are assumed to be smooth.) For connected 
spin manifolds,
the Lichnerowicz method gives a powerful obstruction to this question,
the index $\alpha(M)$ which lives in the (real) K-theory of the
reduced real $C^*$-algebra of the fundamental group.

This invariant just depends on a certain class in a suitably defined 
bordism group (see below). If the dimension of the spin manifold $M$  is bigger than $4$, Gromov and Lawson
\cite{Gromov-Lawson(1980b)} showed that the
existence of a metric with $\scal>0$ on $M$ also just depends on this 
very bordism class. For
some time it was conjectured that in these dimensions the
vanishing of $\alpha(M)$ is also sufficient for the existence 
of such a metric. This conjecture
is called the Gromov-Lawson-Rosenberg (GLR) 
conjecture resp.~the GLR-conjecture for a given group $\Gamma$ 
if we restrict the question to manifolds with this fundamental group. 

The latter is known to be true for a small class of groups, e.g. for the
trivial group \cite{Stolz(1992)}, cyclic groups of odd order 
\cite[1.3]{Rosenberg(1986a)}, $\integers/2$
\cite[5.3]{Rosenberg-Stolz(1994a)} and
quaternionic groups \cite{Botvinnik-Gilkey-Stolz(1997)}. It is also known 
for the free groups, free abelian groups or the fundamental groups of 
orientable surfaces \cite{Rosenberg-Stolz(1995)}. The latter groups 
all have finite dimensional classifying spaces. The first part of 
this paper expands this list, in particular dealing with groups 
with a low dimensional classifying space, and products of such 
groups with free abelian groups.

More precisely:
We show that the statement above is true for $\dim(M)\ge\max\{\dim(B\Gamma)-4,5\}$ if the assembly map $A:KO_*(B\Gamma)\to KO_*(C^*_{red}\Gamma)$
is injective.
In particular the GLR-conjecture is
true in full generality for hyperbolic groups of dimension at most
$9$.
 In addition we show that if $\dim B\Gamma\le 5$ and
the assembly map is split injective for $\Gamma$ then the
GLR-conjecture holds for $\Gamma\times\Z^k$  $\forall k\in\N$.

The invariant $\alpha$ is defined as follows. Let $\pi$ be the 
fundamental group of the given closed connected spin manifold $M$. Then $\alpha(M)$ is the 
image of the bordism class of the classifying map $u:M \to B\pi$ under 
a homomorphism $\alpha:\Omega^{spin}_n(B\pi) \to KO_n(C^*_{red}\pi)$, where
$KO_n(C^*_{red}\pi)$ stands for the $n$-th real $K$-theory group of the 
 real reduced $C^{*}$-algebra $C^*_{red}\pi$, a certain completion of 
the group ring $\R\pi$. The homomorphism $\alpha$ is defined 
through the following factorization
 \begin{equation*}
              \Omega_n^{spin}(B\pi)\stackrel{D}{\to}
ko_n(B\pi)\stackrel{p}{\to} KO_n(B\pi) \stackrel{A}{\to} KO_n(C^*_{red}\pi),
              \end{equation*}
where $ko_*$ is connective real $K$-homology, $KO_*$ the periodic real
$K$-homology, $D$ is the $ko$-theoretic orientation, $p$ the canonical
map between the connective and the periodic theory, and $A$ the
assembly map in topological $K$-theory. Note that for torsion free
groups, this is the Baum-Connes map, and the Baum-Connes conjecture 
states that this map is an isomorphism.
The homomorphism $\alpha$ also has an analytical interpretation 
(cf. \cite[\S\S\ 2--3]{Rosenberg(1986a)}). The image of a class 
$[N \to B\pi]$ in  $KO_{*}(C^*_{red}\pi)$ corresponds to an index 
of an elliptic differential operator on the manifold $N$. The 
Lichnerowicz argument shows that this index must be 
trivial if $N$ admits a metric with positive scalar curvature.
I.~e.~the subgroup $\Omega_n^{spin,+}(B\pi)$ whose elements are 
bordism classes that can be represented by singular manifolds 
with positive scalar curvature is contained 
in the kernel of $\alpha$. However, the counterexample to the 
GLR-conjecture given by Schick \cite{Schick(1998e)} shows that 
the kernel of $\alpha$ is bigger than $\Omega_n^{spin,+}(B\pi)$ 
in general. A weaker version of the GLR-conjecture would ask if at 
least vanishing of $pD[M\xrightarrow{u} B\pi]$ would imply that 
$M$ admits a metric with $\scal>0$.
In the second part of this paper we show that this weakened version 
of the GLR-conjecture fails as well.

If $M$ is not spin, it does not give rise to an element in
$\Omega_{*}^{spin}(B\pi)$ and we cannot apply the above theory.
However, there is a twisted version of all this which can be
applied whenever the universal covering of $M$ is spin. The twist is 
essentially decoded by two cohomology classes $u_1\in H^1(B\pi,\integers/2)$
and $u_2\in H^2(B\pi,\integers/2)$. These classes have to pull back to
the first and second Stiefel-Whitney class of $M$, respectively. 
We then get an obstruction to positive scalar curvature
which lives in $KO_n(C^*_{red}(\gamma(M))$, the $n$-th real $K$-theory 
group of a certain $C^*$-algebra $C^*_{red}(\gamma(M))$ which is 
(up to isomorphism) determined by the triple $(\pi,u_{1},u_{2})$. 
In the third part of the paper we show that the twisted version of the 
(unstable) GLR-conjecture fails as badly as the untwisted one. The example we give  
is a non-orientable manifold which does not admit a metric with $\scal>0$ 
but whose orientation cover actually does have such a metric. 
To show that our manifold indeed is a
counterexample, we use a twisted version of the minimal surface method 
of Schoen and Yau, which relies on geometric measure theory with twists.

\section{The  Gromov-Lawson-Rosenberg conjecture for finite
  dimensional groups}
\begin{theorem}\label{finiteBG}
  Suppose $\Gamma$ is a discrete group and $B\Gamma$ is
  $N$-dimensional. Suppose $n\ge \max({N-4},5)$ and the assembly map
  \begin{equation*}
    A: KO_n(B\Gamma)\to KO_n(C^*_{red}\Gamma)
  \end{equation*}
  is injective. Let $M^n$ be an $n$-dimensional connected closed spin manifold
  with $\pi_1(M)\cong\Gamma$. Then
$M$ admits a metric with positive scalar curvature if and only if its
  analytical index vanishes, i.e.~$\alpha(M)=ApD(M)=0\in KO_n(C^*_{red}\Gamma)$.
\end{theorem}
\forget{\begin{remark}
  For torsion free groups, injectivity of $A$ is  part of the Baum-Connes
  conjecture.\footnote{This part is also called the $KO$-theoretic (strong) Novikov conjecture.} It is known for
  large classes of discrete groups, e.g.~for amenable groups, word-hyperbolic groups, groups of
  finite asymptotic dimension and all discrete subgroups of
  Lie groups (for all this compare
  Higson \cite{Higson(1999)}).
  Higson formulates his results for complex $K$-theory. However, the
  full Baum-Connes conjecture (isomorphism of the Baum-Connes map) in
  the complex case implies isomorphism in the real case. Higson uses
  as basic ingredient of his proof that the full (complex) Baum-Connes
  conjecture is
  true for amenable groupoids. The above principle implies the real version of this
  result. Starting with this and using the real version
  of all the methods Higson uses, one arrives at injectivity of the
  real Baum-Connes map, which we require. (The authors thank Nigel
  Higson for this explanation).
\end{remark}
}%\end{forget}
\begin{remark}\label{grouplist}
  For torsion free groups, injectivity of $A$ is  part of the Baum-Connes
  conjecture.\footnote{This part is also called the $KO$-theoretic (strong) Novikov conjecture.} It is known for
  large classes of discrete groups, e.g.~for amenable groups, 
  word-hyperbolic groups, groups of finite asymptotic dimension 
  and all discrete subgroups of Lie groups (cf.~Higson \cite{Higson(1999)}).
  Higson formulates this results for 
  complex $K$-theory. He uses
  as basic ingredient of his proof that the (complex) Baum-Connes
  conjecture\footnote{The Baum-Connes conjecture states 
  that the Baum-Connes map (which is the assembly map for torsion free groups) 
  is an isomorphism}  is
  true for amenable groupoids. The
   Baum-Connes conjecture in the complex case implies the Baum-Connes conjecture 
  \renewcommand{\thefootnote}{\fnsymbol{footnote}} in 
  the real case. Starting with the real Baum-Connes conjecture for
  amenable groupoids, the descent argument which Higson
  uses implies the strong Novikov conjecture also in the real case,
  which is what we need.\footnote[3]{The authors thank Nigel Higson for explaining 
  the argument to them.}\renewcommand{\thefootnote}{arabic{footnote}}
\end{remark}

\begin{convention}
 If for fixed $\Gamma$ the conclusion of Theorem \ref{finiteBG} 
 is true for any $M$ as above with $\dim M \ge 5$ then we say that the 
 Gromov-Lawson-Rosenberg conjecture is true for $\Gamma$.
\end{convention}

The proof of the theorem above relies on the following theorem of Stephan
Stolz and Rainer Jung
\cite[3.7]{Rosenberg-Stolz(1995)}
\begin{proposition}\label{Stolz}
  Let $M$ be a connected closed spin manifold with 
$\pi_1(M)\cong\Gamma$ and $\dim(M) \ge 5$. If  $D(M)\in ko_*(B\Gamma)$ vanishes, 
then $M$ admits a metric with positive scalar curvature. 
\end{proposition}
\forget{\begin{proof}
Since $0\in ko^+(B\Gamma)$, Stolz main result implies
$M\in \Omega^+(B\Gamma)$.
\end{proof}
} % \end{forget}

To apply this, we first point out some relations between 
spectra and their connective covers. These can easily be verified. 
For the reader's convenience we give the proofs here.

\begin{lemma}\label{compSS}
Let $E$ be a spectrum, $f:E\langle k\rangle\to E$ the $k$-connected cover 
of $E$ ($k\in\Z$) and let $X$ be a CW-complex. Then the induced natural 
transformation  of spectral sequences for the
Atiyah-Hirzebruch spectral sequence
\begin{equation*}
  f^n_{p,q}:E^n_{p,q}(E\langle k\rangle(X))\to E^n_{p,q}(E(X))
\end{equation*}
is an isomorphism for $q\ge k+n-2$ and surjective for $q\ge k$.
\end{lemma}
\begin{proof}
  We prove this by induction on $n$. Since
$  E^2_{p,q}(E(X))=H_p(X,E_q(*))$ and from the definition of $E\langle
k\rangle$ the assertion is true for $n=2$. Now suppose the statement 
is true for some $n\ge 2$. We recall that the $n$-th differentials are a
collection of homomorphisms 
$d^{n}_{p,q}: E^{n}_{p,q}\to E^{n}_{p-n,q+n-1}$, 
and let us denote the differentials in the two spectral sequences by 
$d^{n}_{p,q}\langle k\rangle$ resp.~$d^{n}_{p,q}$.\\
If $q\ge k$ then by assumption
$f^{n}_{p-n,q+n-1}$ is injective and $f^{n}_{p,q}$ is onto. 
Naturality then implies that 
$f^{n}_{p,q}: \Ker d^{n}_{p,q}\langle k\rangle \to \Ker d^{n}_{p,q}$ is
onto, and therefore the same is true for $f^{n+1}_{p,q}$. \\
If $q\ge k+n-1$ then by assumption
$f^{n}_{p+n,q-n+1}$ is onto, and naturality implies that
$f^{n}_{p,q}$ restricts to a surjective map 
$f^{n}_{p,q}: \Image d^{n}_{p+n,q-n+1}\langle k\rangle \to \Image d^{n}_{p+n,q-n+1}$.
Then injectivity of $f^{n+1}_{p,q}$ follows from the injectivity of 
$f^{n}_{p,q}$, which is given by assumption (since $q\ge k+n-2$).
\end{proof}

\begin{lemma}\label{comp_hom}
  Let $E$ be a spectrum, $f:E\langle k\rangle\to E$ the $k$-connected cover 
of $E$ ($k\in\Z$) and let $X$ be a $N$-dimensional CW-complex.
Then the induced natural transformation on the limit terms of the
Atiyah-Hirzebruch spectral sequence
\begin{equation*}
  f^\infty_{p,q}:E^\infty_{p,q}(E\langle k\rangle(X))\to E^\infty_{p,q}(E(X))
\end{equation*}
is  injective for $p+q\ge k+N-1$ and is an isomorphism for the following indices:
\begin{itemize}
\item $q\ge k+N-2$
\item  $p+q\ge k+N-1$ and $q\ge k$
\item $p\ge N+1$ .
\end{itemize}

Consequently, the natural transformation
\[ f_n: E\langle k\rangle_n(X)\to E_n(X) \]
is an isomorphism for $n\ge N+k$ and is injective for $n= N+k-1$.
\end{lemma}
\begin{proof}
  The method of proof is identical to the proof of
  Lemma \ref{compSS},
 proving a corresponding statement for the stages $E^k$,
  and
taking into account that $E^n_{p,q}=0$ for $p>N$. One has to note that
  differentials which map to $E_{0,q}$ are always zero in the
  Atiyah-Hirzebruch spectral sequence. Because
$E^\infty=E^{N}$
the statement for the limit term follows.

If $n\ge N+k$ all terms in the filtration of $E\langle k\rangle_n(X)$ are mapped
isomorphically onto the corresponding terms in the filtration of
$E_n(X)$. Naturality implies that the induced maps on the homology
groups are isomorphisms. Injectivity for $n=N+k-1$ follows from Lemma \ref{Einfty_inj}.
\end{proof}

We want to apply this to $KO$-theory. By definition $ko=KO\langle 0\rangle$. 
Since $\pi_{-1}(KO)=\pi_{-2}(KO)=\pi_{-3}(KO)=0$ we have $ko=KO\langle-3\rangle$.
Hence Lemma \ref{comp_hom} implies 
\begin{proposition}\label{compGhom}
  Suppose $\Gamma$ is a discrete group and  $B\Gamma$ is homotopy
  equivalent to an
  $N$-dimensional $CW$-complex. Then the natural map
  \begin{equation*}
    p_n:ko_n(B\Gamma)\to KO_n(B\Gamma)
  \end{equation*}
  is an isomorphism for $n\ge N-3$ and injective for $n=N-4$.
\end{proposition}

\begin{proof}[Proof of Theorem \ref{finiteBG}]
We only have to show that $p$ and $D$ are injective by
\cite[3.8]{Rosenberg-Stolz(1995)}. We repeat the argument here.
  The ``only if'' part is due to Rosenberg
  \cite{Rosenberg(1986a)}. Suppose now that $ApD(M)=0$. From the assumption
  about the Novikov conjecture it follows $pD(M)=0$. By
  Proposition \ref{compGhom} $p$ is an isomorphism resp.~an injection
  in the dimensions in question. Therefore $D(M)=0$. By the Stolz/Jung result
  \ref{Stolz} the assertion follows.
\end{proof}

\begin{corollary}
  Let $\Gamma$ be a discrete group with $\dim B\Gamma\le 9$. Suppose
  that the
  injectivity part of the Baum-Connes conjecture holds for $\Gamma$. Then the
  Gromov-Lawson-Rosenberg  conjecture is true for
  $\Gamma$. In particular it holds for free groups and surface groups
  and for fundamental groups of complete non-positively curved
  manifolds of dimension $\le 9$.
\end{corollary}
\begin{proof}
  This is true since the GLR-conjecture is just the assertion of
  Theorem \ref{finiteBG} for all $n\ge 5$. For the groups in the list
  the Novikov conjecture is known to be true
  by Remark {grouplist}.
\end{proof}

\forget{\begin{remark}
  In unpublished work, Davis and L\"uck use a spectral sequence
  technique to handle certain groups with torsion. They need very good
  information about the GLR-conjecture for the torsion subgroups, and
  in addition use the same argument we have exploited above conclude.
\end{remark}
}%\end{forget}
\begin{remark}
  To conclude the GLR-conjecture for the free groups or for the 
fundamental groups 
of oriented surfaces we do not need to consider the Atiyah-Hirzebruch 
spectral sequence. Here the injectivity of $p$ follows directly from the 
fact that in these cases $B\pi$ is stably homotopy equivalent to a wedge 
of spheres and from the suspension isomorphism for homology theories 
(cf. \cite[\S 3]{Rosenberg-Stolz(1995)}). This argument can also be 
used to show the GLR-conjecture for all free abelian groups.
\end{remark}

\begin{lemma}\label{timesZ}
If $h_{*}$ is a homology theory then we have  canonical isomorphisms 
$h_{n}(S^{1} \times X) \cong h_{n}(X) \oplus h_{n-1}(X)$ for spaces $X$ 
and $n\in \Z$. In particular, if $f:E\to F$ is a map of spectra then 
$f_*:E_{n}(X\times S^1)\to F_{n}(X\times S^1)$ is injective if and only 
if $f_*:E_{n}(X)\to F_{n}(X)$ and
  $f_{*}:E_{n-1}(X)\to F_{n-1}(X)$ is  injective.\hfill $\Box$ 
\end{lemma}
\forget{ \begin{proof}
  The Leray-Serre spectral sequence of the fibration $X\to X\times
  S^1\to S^1$ shows that we have exact sequences
  \begin{equation*}
    0\to E_n(X)\to E_n(X\times S^1)\to E_{n-1}(X)\to 0
  \end{equation*}
  with a split given by the projection $X\times S^1\to X$. All of
  this, including the split, commute with the natural transformation
  $f_*$. This concludes the proof.
\end{proof}
}%\end{forget}

\begin{proposition}\label{Zerw}
  Suppose $G$ is a group with $\dim BG\le 5$ and the assembly map is
  injective  for $G\times \Z^k$. Then the 
  Gromov-Lawson-Rosenberg conjecture holds for $G\times \Z^k$.
\end{proposition}
\begin{remark}
  If the assembly map for $G$ is a split injection then the same is
  true for $G\times \integers^k$. This can be
  seen using K\"unneth theorems as in Rosenberg's proof of  
  \cite[2.9]{Rosenberg(1986a)} and the fact that the assembly map is 
  an isomorphism for $\Z^k$ \cite[2.6]{Rosenberg(1986a)}. For example, 
  this assumption on $G$ is fulfilled for the groups listed in Remark \ref{grouplist}.
\end{remark}
\begin{proof}[Proof of Proposition \ref{Zerw}]
We show that $ko_*(BG)\to KO_*(BG)$ is an injection for all $*\in\Z$. Then by Lemma \ref{timesZ}
 the same is
true for $ko_*(BG\times \Z^n)\to KO_*(BG\times\Z^n)$ and the proof of
Theorem \ref{finiteBG} applies.

Again, we study the Atiyah-Hirzebruch spectral sequence for $ko_*(BG)$
and $KO_*(BG)$. Since $\dim BG\le 5$ we have
$E^\infty_{p,q}=E^6_{p,q}$. By Lemma \ref{compSS}
$E^\infty_{p,q}(ko(BG))\to E^\infty_{p,q}(KO(BG))$  is an isomorphism for
$q\ge 0$. Moreover, $E^\infty_{p,q}(ko(BG))=0$ for $q<0$.
Now apply the following Lemma \ref{Einfty_inj}.
\forget{If we now solve the extension problems for $ko_n(BG)$ and $KO_n(BG)$
we get (by naturality) isomorphism up to some stage of the
filtration. From this stage on, the quotients are zero for $ko$ and
something for $KO$. Therefore $ko_n(BG)$ remains unchanged and 
therefore is a subgroup of the kernels of the next extensions for $KO_n$.
Eventually $ko_n(BG)$ is mapped injectively into $KO_n(BG)$. 
}%\end{forget}
\end{proof}

\forget{
\begin{lemma}
  \label{Einfty_inj}
  Suppose $\Phi : h_{*} \to h'_{*}$ is a natural transformation of 
  homology theories and for a space $X$ we have two spectral sequences 
  converging to $h_{*}(X)$ and $h'_{*}(X)$ together with a transformation 
  of spectral sequences which is compatible with $\Phi$.
  %(for example, if $X$ is a fibration over a  $CW$-complex) 
  If this transformation is injective on the
  $E^\infty$-term then $\Phi :h_{*}(X)\to h'_{*}(X)$ is injective.
\end{lemma}
\begin{proof}
  To recover $h_{*}(X)$ and $h'_{*}(X)$ from the $E^\infty$-terms, we have 
  to successively solve extension problems 
  \begin{equation*}
    \begin{CD}
          0 @>>> A @>>> B @>>> C @>>> 0 \\
          @VVV    @VV{\Phi}V    @VV{\Phi}V    @VV{\Phi}V   @VVV\\
          0 @>>> A' @>>> B' @>>> C' @>>> 0
        \end{CD}
      \end{equation*}
  where by assumption and induction the maps on kernels and cokernels
  are injective. A diagram chase proves the same for
  $B\to B'$. 
\end{proof}
}%\end{forget}

\begin{lemma}  \label{Einfty_inj}
Let $A$ and $B$ be filtered (graded) groups,
filtered by increasing filtrations $\{A_{s}\}_{s\in\Z}$ and 
$ \{B_{s}\}_{s\in\Z}$. Assume that the corresponding spectral sequence
converges to $A$ (i.e.~$\bigcap_{s\in\integers} A_s =\{0\}$). If a homomorphism $\varphi: A \to B$ respects the filtration and induces injections on the $E^{\infty}$-term, then $\varphi$ is injective.
\end{lemma}

\begin{proof}
For any $s\in \Z$ we have a transformation of short exact sequences
 \begin{equation*}
    \begin{CD}
      1 @>>> A_{s} @>>> A_{s+1} @>>> A_{s+1}/A_{s} @>>> 1\\
       && @VV{\varphi_{s}}V    @VV{\varphi_{s+1}}V  
         @VV{\overline{\varphi}_{s+1}}V\\
      1 @>>>  B_{s} @>>> B_{s+1} @>>> B_{s+1}/B_{s} @>>> 1.\\
        \end{CD}
      \end{equation*}
$\overline{\varphi}_{s+1}$ is one of the injections in the $E^{\infty}$-term, hence we have 
$\Ker \varphi_{s} = \Ker \varphi_{s+1}$. Therefore
$$\Ker \varphi_{s} = \bigcap_{s' \in \Z} \Ker \varphi_{s'} \subset 
     \bigcap_{s' \in \Z} A_{s'} = 0.$$
I.e.~all the $\varphi_{s}$ are injective and therefore so is 
$\varphi = \lim_{s\to \infty} \varphi_{s}$.
\end{proof}

\forget{
\begin{remark}
  Proposition \ref{Zerw} in particular  recovers the result of Rosenberg and Stolz
  \cite{Rosenberg-Stolz(1995)} that the GLR-conjecture is true for
  free abelian fundamental groups, and  more generally
  for products of free abelian groups with fundamental groups of complete
  non-positively curved manifolds of dimension $\le 4$.
\end{remark}
}%\end{forget}

% end{forget}

\begin{proposition}
  The GLR-conjecture is true for every group $G$ which fits into an
  exact sequence $1\to\integers^n\to G\to \integers\to 1$. 
\end{proposition}
\begin{proof}
  We study the Atiyah-Hirzebruch spectral sequence for the
  corresponding fibration of classifying spaces $T^n\to BG\to S^1$. In the case of a
  fibration over $S^1$, we have an action of $\integers$ on the
  homology of the fiber (induced by the homotopy class of self-maps of
  the fiber which comes from fiber transport along a generator of
  $\pi_1(S^1)$),  and the $E^2$-term is given by the invariant
  and coinvariant elements, respectively.

  Note that we can compute $ko_p(T^n)$ and $KO_p(T^n)$ using the
  K\"unneth theorem, and that each of these groups is a product of
  subgroups which have a given ``suspension''-dimension, where an
  element in $ko_{p-k}(*)\times [S^1]_{i_1}\times \dots [S^1]_{i_k}$
  has suspension-dimension $k$. Here $[S^1]_{i}$ is a generator of
  the first homology of the $i$-th factor of $T^n$.  Since fiber transport
  gives rise to a map of the first
  homology of $T^n$, the induced map in $K$-theory will preserve the
  suspension dimension.

  The transformation $ko_p(T^n)\to KO_p(T^n)$ is exactly the inclusion
  of those summands with suspension dimension less than $n-p+1$. 

We get the following $E^2$-terms:
  \begin{equation*}
    \begin{split}
      E^2_{0,q}(ko)= ko_q(T^n)/\integers &\xrightarrow{p}
      KO_q(T^n)/\integers=E^2_{0,q}(KO)\\
      E^2_{1,q}(ko)=
    ko_q(T^n)^\integers &\xrightarrow{p} KO_q(T^n)^\integers = E^2_{1,q}(KO).
    \end{split}
  \end{equation*}
  All the other groups are trivial. In particular, 
  the $E^2$-terms are already the $E^\infty$-terms. 
  Since $ko_*(T^n)$ is a $\integers$-invariant
  summand in $KO_*(T^n)$, the same is true for the groups of invariant
  and coinvariant elements, respectively. In particular, the natural
  transformation induces a (split) injection of the $E^\infty$-terms
  for $ko$ and $KO$. Lemma \ref{Einfty_inj} concludes the proof.
\end{proof}

\begin{corollary}
  If $G$ is the fundamental group of a closed Riemannian manifold
  whose holonomy is  $\integers/p$ for a prime $p$ then the
  GLR-conjecture is true for $G$.
\end{corollary}
\begin{proof} 
Such a group $G$ fits into an exact sequence as in the above proposition 
(cf. \cite[{\small $\binom{*}{*}$} on page 552]{Charlap-Vasquez(1965)}).
\end{proof}

\forget{
\begin{remark}
The results of this section give certain implications of the
injectivity part of the Baum-Connes conjecture. Powerful obstructions
to positive scalar curvature therefore could in principle lead to
counterexamples to the injectivity of the Baum-Connes map. However,
the difficult side always seems to be the $K$-theory of the
$C^*$-algebra, so we don't gain very much. Moreover, so far the only
(additional) obstruction to $\scal>0$ which is known is the minimal
surface method combined with the Gauss-Bonnet theorem, which means we
have to find an example where injectivity of the Baum-Connes map fails
in dimension $2$. But the $2$-dimensional case is by work of Valette
\Kommentar{da gibt es einige Koautoren --- heraussuchen!}  quite well
understood and it's not unlikely that one can prove that
$2$-dimensional classes are always mapped injectively by the
Baum-Connes map.
\end{remark}
}% end{forget}

\section{A counterexample to a weak form of the
  Gromov-Lawson-Rosenberg conjecture}

In this section we present a closed connected spin manifold $M$ with
classifying map
$M\xrightarrow{u}BG$ such that $pD[M\xrightarrow{u}BG]=0$, but $M$ does
not admit a metric with positive scalar curvature. To do this, we
start with a closer look at the $KO$-homology of finite groups.

\begin{proposition}\label{KO}
  If $G$ is a finite group of odd order then $\widetilde{KO}_2(BG)=0$.
\end{proposition}
\begin{proof}
Note the following two facts:
\begin{itemize}
\item $\widetilde{KU}_2(BG)=0$ by \forget{\cite[p. 615]{Vick(1970)} or }
\cite[7.1]{Greenlees(1993a)}.
\item Complexification and forgetting the complex structure yields two
natural transformation $c: KO\to KU$ and $r:KU\to KO$ whose
composition $rc: KO\to KU\to KO$ is
multiplication by $2$.
\end{itemize}
These two facts immediately imply that on $\widetilde{KO}_2(BG)$ multiplication
by two is equal to zero, i.e. $\widetilde{KO}_2(BG)$ is a $\Z/2$-vector space.
\forget{
\begin{proposition}
If $G$ is a finite group of odd order $n$ then $\widetilde H_*(G,\Z)$  is
$n$- torsion.
In particular, multiplication by two is an isomorphism on $\widetilde H_*(G,\Z)$.
\end{proposition}
\begin{proof} For convenience, we repeat the proof:\\
  The contractible space $EG$ is an $n$-sheeted covering of
  $BG$. Therefore multiplication by $n$ on $\widetilde H_*(BG,\Z)$, which is the
  composition of the transfer map to $\widetilde H_*(EG,\Z)$ with the
  projection map $\widetilde H_*(EG)\to \widetilde H_*(BG)$ factorizes through
  the trivial group $\widetilde H_*(EG)$. Therefore, it is zero.
\end{proof}
} %end{forget}
On the other hand it is well-known that  $\widetilde H_*(G,\Z)$ is annihilated 
by the order of $G$. Hence  $\widetilde H_*(G,\Z)$ is odd torsion and 
multiplication by $2$ is an isomorphism on  $\widetilde H_*(G,\Z)$.
Using the Atiyah-Hirzebruch spectral sequence and the $5$-lemma, we see
that multiplication by 2 is an isomorphism on the $\Z/2$-vector
space  $\widetilde{KO}_2(BG)$. Therefore this group must be
trivial.
\end{proof}

\begin{proposition}
  Let $G$ be a discrete group. Then
  \begin{equation}
    \label{omega2plus}
    \Omega^{spin,+}_2(BG):=\left\{
       \text{\begin{tabular}{c}
              bordism classes $[M\to BG] \in\Omega^{spin}_2(BG)$,\\
              where $M$ admits a metric with $\scal>0$
             \end{tabular}} 
                                 \right\}= {0} .
  \end{equation}
\end{proposition}
\begin{proof}
 By the Gauss-Bonnet theorem there is only one orientable $2$-manifold with positive
(scalar) curvature, namely $S^2$. Since $\pi_2(BG)$ is trivial, up to homotopy
only the trivial map from $S^2$ to $BG$ exists. Therefore only the trivial
element in $\Omega_2^{spin}(BG)$ can be represented by a manifold with
positive scalar curvature.
\end{proof}

\begin{proposition}\label{kernel_not_zero}
  For $G=\Z/3\times \Z/3$, the index map
  \begin{equation}
    \label{ker_ind}
    pD: \Omega^{spin}_2(BG)\to KO_2(BG)
  \end{equation}
has a non-trivial kernel.

In particular, not every element of this kernel lies in $\Omega^{spin,+}_2(BG)$.
\end{proposition}
\begin{proof}
It is sufficient to show this for the reduced groups. Since by Proposition
\ref{KO} $\widetilde{KO}_2(BG)=0$, we only have to show that
$\widetilde\Omega^{spin}_2(BG)$ is non trivial. This follows immediately from the
Atiyah-Hirzebruch spectral sequence since 
$E^2_{2,0}\cong H_2(B\Z/3\times B\Z/3,\Z)\ne 0$, and the only possible 
differential at $(2,0)$ is the $E^2$-differential mapping to $E^2_{0,1}$.
This differential is zero because the coefficients survive to $E^\infty$.
\end{proof}

\begin{theorem}\label{cex}
  There is a $5$-dimensional closed connected spin manifold $M$ with
  $ \pi_1(M):=\pi\cong \Z^3\times\Z/3\times\Z/3$ such that
  \begin{equation*}
    pD[M]=0 \in KO_5(B\pi),
  \end{equation*}
but $M$ does not admit a metric with positive scalar curvature.
\end{theorem}
The proof, using the minimal surface technique of Schoen and Yau,
 will be given in the next section.
\begin{remark}
  This is a counterexample to one direction of an older conjecture
  of Gromov and Lawson\forget{\cite{Gromov-Lawson(1983)}}, 
  namely that positive
  scalar curvature is equivalent to the vanishing of $pD[M]$. The
  other direction of the conjecture was disproved by Rosenberg 
  \cite{Rosenberg(1983)}. He then suggested
  that vanishing of the image $\alpha(M)$ of $pD[M]$ under the assembly 
  map might be sufficient for the existence of a metric with 
  $\scal>0$ on $M$. In the following this has become known 
  as the GLR-conjecture.

Dwyer and Stolz have another counterexample like $M$ with the
additional very nice property
that $\pi_1(M)$ is torsion-free. Their computation is based on
homotopy theoretic considerations and the Baumslag-Dyer-Heller
construction of groups with given homology
\cite{Baumslag-Dyer-Heller(1980)}. The drawback of this approach is
that the group they construct is rather artificial.
\forget{ (The original
version of the GLR-conjecture requires  that the groups in
question are ``reasonable'', and most definitions of ``reasonable''
probably exclude the group of Dwyer and Stolz.)
}% end{forget}
\end{remark}

\section{The minimal surface obstruction to positive\\ scalar curvature}

The following
theorem is the differential geometrical backbone for the application of
minimal surfaces to the positive scalar curvature problem:
\begin{theorem}\label{diffgeo}
Let $(M^m, g)$ be a Riemannian manifold with $\scal>0$, $\dim M=m\ge 3$.
 If $V$ is a smooth $(m-1)$-dimensional
immersed submanifold of $M$ with trivial normal bundle,
 and if $V$ is a local minimum of the
$(m-1)$-volume, then $V$ admits a metric with $\scal>0$, too.
Note that orientability is not required.
\end{theorem}
\begin{proof} \cite[Proof of Theorem 5.1]{Schoen-Yau(1979a)} for $m=3$, 
and \cite[Proof of Theorem 1]{Schoen-Yau(1979)} for $m>3$. The
computations there are formulated only for $M$ and $V$ oriented, but
they carry over literally to the case where $V$ has trivial normal
bundle in $M$.
\end{proof}

%{Das folgende Theorem ist das was wir wirklich brauchen.}
\begin{theorem}\label{geomeasure}
  Suppose $M$ is a closed Riemannian manifold with $\dim M=m\le
  8$. Furthermore let $\alpha\in H^1(M,\Z)$, and let $\obundle$ denote
the orientation bundle of $M$. Then
  \[ x:=\alpha\cap [M]\in H_{m-1}(M,{\obundle}) \]
can be represented by an embedded hypersurface $V$ with trivial 
normal bundle which is a local minimum for
$(m-1)$-volume (if $m=8$ with respect to suitable metrics arbitrarily
close in $C^3$ to the metric we started with).

 Note that $V$ also represents $\overline{x} \in H_{m-1}(M,\integers/2)$, the $mod\ 2$ reduction of $x \in H_{m-1}(M,{\obundle})$.
\end{theorem}
\begin{proof}
If $M$ is orientable then we do not have to twist the coefficients and
for $m\le 7$
this is a classical result of geometric measure
theory (cf. \cite[Chapter 8]{Morgan(1988)}). The case $m=8$ follows 
from the following result of Smale \cite{Smale(1993)}: the
set of $C^k$-metrics for which the regularity statement holds is open and dense in the set
of all $C^k$-metrics ($k\ge3$ and with the usual Banach-space topology). We are
only interested in
$C^\infty$-metrics. But these are dense in the set of $C^k$-metrics,
which concludes the proof for the orientable case.

More generally, if we work with homology and cohomology with twisted
coefficients in the orientation bundle and its dual, respectively,
then every submanifold with orientable normal bundle represents a
homology class in this theory (being dual to the corresponding
differential forms), and there is a Poincar\'e\ duality between this
homology theory and untwisted integral cohomology. Now one can
tediously check that all the basic results of geometric measure theory
are true for this theory with twisted coefficients, including the work
of Smale. In particular, for
each twisted $(m-1)$-cohomology class we find a representing current
with minimal volume, and such a minimizing current is a smooth
embedded hypersurface with trivial normal bundle if $m\le 7$, and we
can do this for a generic metric if $m=8$.
Details of this will be presented in a forthcoming
paper.
\end{proof}

Remember that if we are given a class $\alpha \in H^{1}(M,\Z)$ we may 
represent it by a map $f:M \to S^{1}$ being transverse to $1\in S^{1}$.
Then $V=f^{-1}(1) \subset M$  represents $\alpha \cap [M]$ (and
conversely, every hypersurface representing $\alpha \cap [M]$ is
obtained in this way). Furthermore,
if $f':M\to S^{1}$ is a second map as above and $V' = {f'}^{-1}(1)$ then 
$f$ and $f'$ are homotopic, and a homotopy $H: f\simeq f'$ being transverse to 
$1\in S^{1}$ provides a bordism $W= H^{-1}(1): V \sim V'$ embedded in
$M\times [0,1]$. If $\Omega$ denotes some $B$-bordism theory%
\footnote{for example Spin bordism, oriented bordism or unoriented
  bordism}\forget{(cf.~\cite{Stong(1968)})} and $M$ is equipped with a (normal)
$B$-structure then the above procedure provides a class in
$\Omega_{m-1}(M)$. This follows from the fact that the normal bundle
of $V \subset M$ respectively $W \subset M\times [0,1]$ is trivial.

More generally, let $X$ be a space and $\Phi: W \to X$ be a $B$-bordism 
between singular $B$-manifolds $\phi: M \to X$ and $\phi': M' \to X$. 
If $f: X \to S^{1}$ represents an element $\alpha \in H^{1}(X,\Z)$, 
then $f \circ \Phi$ is homotopic to a map $\Psi: W \to S^{1}$ with
$\Psi$ and $\Psi|\partial W$ being transverse to $1\in S^{1}$
(moreover, the map $\Psi|\partial W$ with the corresponding properties
may be given in advance). 
Then $\Psi^{-1}(1) \subset W$ is a $B$-bordism between the 
hypersurfaces $(\Psi)^{-1}(1) \cap M \subset M$ and 
$\Psi^{-1}(1) \cap M'\subset M'$, and restricting $\Phi$ to  
$\Psi^{-1}(1)$ now yields a singular $B$-bordism between 
singular hypersurfaces into $X$. \\
If $f':X\to S^1$ is homotopic to $f$,  a similar construction as above gives a
singular $B$-bordism between the resulting singular hypersurfaces into $X$.

 Taken together, we obtain a homomorphism
\begin{equation}
  \label{pairing}
 \cap: H^1(X,\integers)\times\Omega_m(X) \to \Omega_{m-1}(X) .
\end{equation}

The two theorems above now imply
\begin{theorem}\label{bord_schoen_yau}
  Let $X$ be a space and let $3\le m\le 8$. Then \eqref{pairing} 
restricts to a homomorphism:
  \begin{equation}
    \label{pairing+}
     H^1(X,\integers)\cap  \Omega_m^+(X)\to \Omega_{m-1}^+(X)
  \end{equation}
where $ \Omega_*^+(X) \subset  \Omega_*(X)$ is the subgroup of 
bordism classes which can be represented by singular manifolds 
which admit a metric with $\scal >0$.
\end{theorem}
\begin{proof}
  Theorem \ref{geomeasure} implies that, given $f:M\to S^1$ (dual to a
  given class in $H_{m-1}(M,\obundle)$) we find a
  homotopic map $g:M\to S^1$ which is transverse to $1$ and such that
  the hypersurface $V=g^{-1}(1)$ is minimal for the $(m-1)$-volume (in
  dimension $8$ we  replace the given metric
  by one which
  is $C^3$-close). In any case, since the scalar curvature is
  continuous with respect to the $C^3$-norm on the space of all
  Riemannian metrics, $V$ is volume minimizing with respect to a metric with
  positive scalar curvature whenever we start with such a metric. By
  Theorem \ref{diffgeo}, it admits a metric with $\scal>0$.
\end{proof}

The spin version of \eqref{pairing+} is used in Schick
\cite{Schick(1998e)} to give a counterexample to the
GLR-conjecture. We are now ready to use \eqref{pairing+} to complete
the proof of the theorem in the previous section.

\begin{proof}[Proof of Theorem \ref{cex}]
  We use the same construction as in \cite{Schick(1998e)}: take a two 
  dimensional singular
  bordism $f: F\to B(\Z/3)^2$ which does not lie in
  $\Omega_2^{spin,+}(B(\Z/3)^2)$ with the property that the image of its
  $KO$-orientation is zero in $KO_2(B(\Z/3)^2)$.
  We have already shown that such an $f$ exists
  (see \ref{kernel_not_zero}). Then take the
  product with $id: S^1\times S^1\times S^1\to
  B\Z^3$. The image of the $KO$-orientation of the new singular bordism in 
  $KO_5(B\pi)$ is still zero. Doing surgeries we find a spin bordant manifold 
  $M$ with fundamental group $\pi$ such that $[f\times id] = [u:M \to B\pi]$ 
  with $u:M \to B\pi$ being the classifying map of the universal covering. 
  I.e.~$M$ satisfies $\alpha(M) = 0 \in KO_5(B\pi)$.

  Now the Schoen-Yau method   shows that $M$ cannot admit a metric with
  positive scalar curvature because taking successively minimal
  surfaces dual to the three canonical generators of $H^{1}(B\Z^{3},\Z)$ 
  otherwise would yield a manifold bordant to $F$ which would also admit 
  a metric with $\scal>0$. But this contradicts our choice of $F$.
\end{proof}

\section{The example for non-orientable minimal surfaces}\label{sec_ex}

We want to show that \eqref{pairing+} yields interesting results for
non-orientable manifolds. If such a manifold has a metric with
positive scalar curvature then the same is true for its
orientation cover. The converse is not necessarily true.
We will produce a five dimensional example whose orientation cover
admits a metric with $\scal>0$.

We start with the Klein bottle $K$, the
non-orientable surface whose orientation cover is the torus. It follows 
from the Wu relations that $w_{2}(K) = w_{1}(K)^{2}$. Hence the tangent 
bundle $TK$ admits a $Pin(2)$-structure 
\forget{\footnote{$Pin(n) \subset \Cliff_{n}$ 
is generated by the elements of $S^{n-1} \subset \R^{n}\subset \Cliff_{n}$.
 $\Cliff_{n}$ is obtained from the tensor algebra of $\R^n$ 
 subject to the relations $v\cdot v = - |v|^2 \cdot 1, v\in \R^{n}$. }
}%\end{forget}
(cf. \cite{Giambalvo(1973)}).

\begin{lemma}\label{pin_spin} If $W$ is an $n$-dimensional
  $Pin(n)$-manifold then its orientation cover $\widetilde{W}$
  canonically carries a $Spin(n)$-structure.
\end{lemma}

\begin{proof} If $P \to W$ is a $Pin(n)$-principle bundle then 
$P \to P/Spin(n)$ is a $Spin(n)$-principle bundle over $P/Spin(n)$, 
and $q: P/Spin(n) \to W$ is a double covering classified by $w_{1}(W)$, 
i.e.~we may regard $P/Spin(n)$ as a model for $\widetilde{W}$. 
Moreover, if $(P,f: P\times_{Pin(n)}\R^{n} \cong TW)$ represents a 
$Pin(n)$-structure on $W$ then there is a canonical map 
$\widetilde{f}: P\times_{Spin(n)}\R^{n} \cong q^{*}TW = T\widetilde{W}$ covering $f$. Hence $(P,\widetilde{f})$ represents a $Spin(n)$-structure on $\widetilde{W}$.
\end{proof}

It is an easy exercise to construct a $Pin(2)$-structure on $K$ such that $\widetilde{K} = S^1\times S^1$ is the product of two copies of $S^1$ carrying the zero bordant $Spin(1)$-structure. In the following we shall think of $K$ being equipped with this $Pin(2)$-structure.\\

Next recall that
$\pi_1(K)= \Z\ltimes\Z$ is the non-trivial semidirect product of $\Z$
with $\Z$, where the first factor is generated by the standard loop in
$K$ with non-trivial normal bundle, and the second factor (which is the
normal subgroup) by the loop with trivial normal bundle. The image of
the fundamental group of the orientation cover is then $(2\integers)\ltimes\integers=(2\integers)\times\integers$. Let $p=(p_1,p_2): \pi_1(K)\to\Z/2\times \Z/2$ be the projection 
corresponding to the structure of the semidirect product (i.e.~the
image of $p_2$ is generated by $p$ applied to the loop with
trivial normal bundle, and $p_1$ represents the first Stiefel-Whitney class).
 Let us use the same letters for the corresponding maps 
\[ p = (p_{1},p_{2}) : K\to B\Z/2 \times B\Z/2 .\]
Recall also that a $Pin(n)$-structure on the tangent bundle of a manifold 
corresponds to a $Pin'$-structure on its stable normal bundle 
(cf.~\cite{Giambalvo(1973)}). 
In particular,
$p_{2}: K \to B\Z/2$ represents a bordism class 
$[p_{2}] \in \Omega^{pin'}_{2}(B\Z/2)$. Since the product of a
$Pin(n)$-manifold and a $Spin(m)$-manifold canonically carries a 
$Pin(n+m)$-structure, we can multiply $[p_{2}]$ with 
$[id: S^1\to S^1]^3\in \Omega_{3}^{spin}(B\Z^3)$, and we 
obtain $[p_{2}\times id] \in \Omega_{5}^{pin'}(B\pi')$ with 
$\pi' = \Z/2\times \Z^3$. $p_{2}\times id$ has a classifying map
$$(\hat{\nu}_{N}, p_{2}\times id): N = K \times \left( S^{1}\right)^{3}
\longrightarrow BPin' \times B\pi'.$$
In $H^{1}(BPin',\Z/2) \cong [BPin',B\Z/2]$, let $w: BPin'\to B\Z/2$
represent the non-trivial element. Then $w\circ\hat{\nu}_{N}$
represents $w_1(N)$,
i.e.~$w\circ\hat{\nu_N}$ is homotopic to $p_1 \circ pr_{K}: N \to K \to B\Z/2$.

Now we can do surgeries on $N \times [0,1]$ to obtain a singular 
$BPin'$-bordism $\Phi: W \to B\pi'$ between $(p_{2}\times id)$ and a 
singular $BPin'$-manifold $\phi: M \to B\pi'$ such that the classifying map
$$(\hat{\nu}_{M},\phi): M \longrightarrow  BPin' \times B\pi'$$
is a $2$-equivalence. Then $w\circ\hat{\nu}_{M}$ represents $w_1(M)$.

Note that since $\phi$ and $p_{2}\times id$ are $pin'$-bordant 
so are
$(\hat{\nu}_{M},\phi)$ and $(\hat{\nu}_{N}, p_{2}\times id)$, and hence so are
$$ u = ( w\circ\hat{\nu}_{M}, \phi) : 
M \longrightarrow B\Z/2 \times B\pi' = B\pi,\quad \pi = \Z/2 \times
\pi'$$
and $(w \circ \hat{\nu}_{N}, p_{2}\times id) \simeq (p_{1}, p_{2} \times id) = (p \times id)$, i.e.~we have $[u] =  [p\times id] \in \Omega_{*}^{pin'}(B\pi)$.

\begin{proposition}\label{proposition_in_section_4}\
\begin{itemize}
\item[{\rm (1)}] The orientation covering $\widetilde{M}$ does admit a 
metric with $\scal>0$.
\item[{\rm (2)}] $M$ does not admit a metric with $\scal>0$.
\end{itemize}
\end{proposition}

\begin{proof}

(1)
Since $(\hat{\nu}_{M},\phi): M \to BPin' \times B\pi'$ is a 
$2$-equivalence it induces an isomorphism $\pi_{1}(M) \cong \Z/2 \times \pi'$. 
Therefore the map $ u : M \longrightarrow  B\pi$
above classifies the universal covering of $M$. Since $w\circ\hat{\nu}_{M}$
represents $w_{1}(M)$ the universal cover of the orientation covering 
$\widetilde{M}$ is classified by $\widetilde{u} = \phi \circ q$, 
where $q: \widetilde{M} \to M$ is the orientation covering. 
Using the bordism results of Gromov-Lawson mentioned in the 
introduction it follows that $\widetilde{M}$ has a metric with $\scal >0$ 
if and only if $[\widetilde{u}]\in \Omega_{5}^{spin,+}(B\pi')$. 
However, using the orientation cover of the bordism $W$ between $M$
and $N$ and Lemma \ref{pin_spin}, $\widetilde{u}$ is spin bordant to 
$$S^{1}\times S^{1} \times (S^1)^3 \to K \times (S^1)^3 
\xrightarrow{p_2\times id} B\Z/2\times B\left( \Z^{3} \right) = B\pi'.$$ 
This is zero bordant since we can fill in a disk $D^2$ into 
the first copy of the $S^1$'s and still get a map to $B\pi'$, 
i.e.~$[\widetilde{u}]$ lies in the subgroup $\Omega_5^{spin,+}(B\pi')$.

(2)
Consider the composition $u: M\to B\pi$ above and let us forget about 
the $Pin(5)$-structure of $M$. Then $u$ represents a bordism class 
$[u]\in\Omega_{5}^{o}(B\pi)$, the unoriented bordism group of 
$B\pi$. Moreover $p\times\id:N\to B\pi$  represents the same
element. We now show that $[u]=[p\times id] \notin \Omega_{5}^{o,+}(B\pi)$.
To do this we cap the $\Z/2$-orientation 
class with elements of $H^1(N,\Z)$ which come from $B\pi'$. 
Doing this three times with the generators of 
$H^1(B\Z^3,\Z) \subset H^1(B\pi',\Z)$ 
eventually yields the bordism class of
\[ K\stackrel{p}{\to} B\Z/2\times B\Z/2\to B\pi .\]
We have to show that this is not in $\Omega_2^{o,+}(B\pi)$.
It suffices to show that $[p:K\to B\Z/2\times B\Z/2]$ is not in
$\Omega_2^{o,+}(B\Z/2\times B\Z/2)$.

 By the Gauss-Bonnet theorem, elements in $\Omega_2^{o,+}(X)$ for a 
space $X$ are represented  by maps either from $S^2$ or $\R P^2$.
Since $K=\R P^2\#\R P^2$ is null bordant, it is not
bordant to $\R P^2$. Hence it remains to rule out that $p$ 
is bordant to some map $f:S^2\to B\Z/2 \times B\Z/2$. However, every map 
$S^2 \to B\Z/2 \times B\Z/2$ is homotopically trivial and 
$S^{2} = \partial D^2$. I.e.~it is enough to show that 
$[p] \neq 0 \in \Omega_2^{o}(B\Z/2\times B\Z/2)$. To do so we compute 
a non-trivial characteristic number of $[p]$.

Let $H^1(B\Z/2 \times B\Z/2,\Z/2) = \Z/2\langle u_1\rangle 
\oplus \Z/2\langle v_1 \rangle$ be induced by the K\"unneth isomorphism. 
The two classes $u_1$ and $v_1$ pull back to two generators $p^*u_{1}$ 
and $p^*v_1$ in $H^1(K,\Z/2)$. The latter are dual to the two standard 
loops generating the fundamental group. These intersect in exactly one point. 
Hence their intersection number is non-trivial, which implies that 
$<p^*u_{1} \cup p^*v_1, [M]> \neq 0$. Hence we have found a 
non-trivial characteristic number of $[p]$, which now completes the
proof, taking Theorem \ref{bord_schoen_yau} into account.
\end{proof}

\section{Twisted index of manifolds with universal covering spin}

The Lichnerowicz method for positive scalar curvature obstructions a priori
works only for spin manifolds. Rosenberg \cite[3.4]{Rosenberg(1986a)}
proved the following powerful version:
\begin{theorem}\label{Rosenbergs.obstruction}
  If $M$ is an $n$-dimensional connected closed spin manifold 
with positive scalar curvature,
  then the index
    $\alpha(M)= 0\in KO_n(C^*_{red}\pi_1(M))$.
\end{theorem}

In \cite{Rosenberg(1983)} Rosenberg explains how one can extend the theory
to manifolds which are not necessarily spin, but where the universal
covering is spin. This was later refined by Stephan Stolz to take
non-orientable manifolds into account. The definite reference for 
these generalizations is \cite{Stolz(1998)}.\\

 Let $M$ be a connected manifold with fundamental group $\pi$. 
If the universal 
covering of $M$ is spin then the first two Stiefel-Whitney classes of 
$M$ pull back from two cohomology classes $u_{1}\in H^1(B\pi,\Z/2)$ and  
$u_{2}\in H^2(B\pi,\Z/2)$ via a map $u:M \to B\pi$ classifying the 
universal covering. $u_{1}$ can be interpreted as a homomorphism 
$\pi \to \Z/2$ and $u_{2}$ determines (up to isomorphism) a central 
extension 
$$1\longrightarrow \Z/2\longrightarrow \hat{\pi}\longrightarrow 
\pi \longrightarrow 1.$$

\begin{definition}
A supergroup $\gamma=(\pi, w,\hat{\pi})$ is a triple consisting 
of a group $\pi$, a homomorphism $w:\pi \to \Z/2$ and a central extension
as above.
\end{definition}

Given a supergroup $\gamma$, Stolz constructs a corresponding bordism theory 
$\Omega_{*}(\gamma)$, a $C^*$-algebra $C^{*}_{red}(\gamma)$ and homomorphism
$\alpha: \Omega_{*}(\gamma) \to KO_{*}(C^{*}_{red}(\gamma))$.
In the case where $w$ is trivial and the extension is split the 
constructions give back the homomorphism 
 $\alpha:\Omega^{spin}_*(B\pi) \to KO_*(C^*_{red}\pi)$ used above. 
The bordism group is given by a (normal) $B$-bordism theory with 
$B=B(\gamma)$ constructed from the supergroup data. Stolz then shows 
that a connected closed manifold $M$ determines a supergroup $\gamma(M)$ such 
that $M$ has a $B(\gamma(M))$-structure whose classifying map 
$M \to B(\gamma(M))$ is a $2$-equivalence
 if $\dim M = n \ge 3$.\footnote{Note that $B(\gamma(M))$
does not coincide with $BG(n,\gamma(M))$ as defined in 
\cite[\S 2]{Stolz(1998)} 
since we describe $\Omega_{*}(\gamma)$ through manifolds with a structure 
on its stable {\it normal} bundle. However, $M \to B(\gamma(M))$ is a 
2-equivalence if and only if $M\to BG(n,\gamma(M))$ is one 
%(cf.~\cite[2.12]{Stolz(1998)}) 
(for $n \ge 3$).}
If $n\ge 5$ then it  follows from the surgery results of 
Gromov-Lawson (cf. \cite[3.3]{Stolz(1995)}) that $M$ has a 
positive scalar curvature metric if and only if 
$[M]:=[id_{M}] \in \Omega_{n}(\gamma(M))^{+} \subset \Omega_{n}(\gamma(M))$, 
which actually motivated the definition of $\Omega_{*}(\gamma)$.\\

On the other hand we have the following generalization 
of Theorem \ref{Rosenbergs.obstruction}:
\begin{theorem}
If $M$ has a metric with $\scal >0$ then $\alpha([M]) = 0$.
\end{theorem}
\begin{proof} This is a consequence of \cite[1.1 and 1.3]{Stolz(1998)}.
As in the spin case  it follows from an index theoretic 
description of the homomorphism $\alpha$ combined with a Lichnerowicz 
type argument.
\end{proof}

We now show that the manifold $M$ defined in the previous section provides
a non-spin counterexample to the generalized GLR-conjecture, which would say
that a closed manifold $M$ of dimension $n\ge 5$ has a metric with $\scal >0$ 
if and only if $\alpha([M]) = 0 \in  KO_n(C^*_{red}(\gamma(M)))$.\\

First we need to work out the associated supergroup $\gamma(M)$. By definition 
$\pi = \pi_{1}(M)$ and $w=w_{1}(M)$. The definition of the extension
follows from

\begin{lemma}
If $w_{2}(M) = u^{*}(e)$ for $u: M\to B\pi$ a map classifying the universal 
covering and $e\in H^2(B\pi,\Z/2)$ then 
$\gamma(M) = (\pi,w_{1}(M), \hat{\pi})$ with $\hat{\pi}$ an 
extension classified by $e$.
\end{lemma}
\begin{proof} \cite[2.11 (2)]{Stolz(1998)}.\end{proof}

In our example we have $w_{2}(M) = u^{*}(u_{1}^2)$ using the notation 
introduced in the preceding section. It follows that
$$\gamma(M) = (\pi = \Z/2\times \pi', pr_{1}:\pi \to \Z/2, Pin(1)\times \pi').
$$
Going through the definitions we find $B(\gamma(M)) = BPin'\times B\pi'$. 
Hence we have an isomorphism 
$\Omega_{*}(\gamma(M)) \cong \Omega_{*}^{pin'}(B\pi')$.\\

On the other hand the corresponding $C^*$-algebra $C^{*}_{red}(\gamma(M))$ is 
given by
\begin{equation}\label{tensor.product} 
C^{*}_{red}(\gamma(M))  \cong \Cliff_{1} \hat{\otimes} C^{*}_{red}(\pi'),
\end{equation}
where the tensor product we use here is the $C^*$-tensor product, which 
is uniquely defined in our case.
To obtain this isomorphism  we may regard $\gamma(M)$ as a product of 
supergroups \cite[2.4]{Stolz(1998)}
$$\gamma(M) = \gamma_{1} \hat{\times} \gamma_{2} \text{\ \  with \ \ }
\begin{array}{l}
 \gamma_{1} = (\Z/2, id, Pin(1))\\
 \gamma_{2} = (\pi', 0, \pi'\times \Z/2) .
\end{array}$$
Then by \cite[8.2 (2) and 8.3]{Stolz(1998)} we obtain
$C^{*}(\gamma_{1}) \cong \Cliff_{1}$ and
 $C^{*}(\gamma_{2}) = C^{*}_{red}(\pi')$,
and \cite[8.4]{Stolz(1998)} then implies \eqref{tensor.product}. Note
that all the way through we deal with $\integers/2$-graded
$C^*$-algebras, where $C^*_{red}(\Gamma)$ for a group $\Gamma$  is
trivially graded and
$\Cliff_1$ has the standard grading.

\begin{lemma} $KO_5(C^*_{red}(\gamma(M))) \cong\Z^6$.
\end{lemma}

\begin{proof}  First note that we deal with abelian groups. 
Therefore reduced and unreduced $C^*$-algebras coincide.
%(see Schroeder: remarks after 1.3.2)
Hence we have
\begin{equation*}\begin{split}
  KO_n(C^*_{red}(\gamma(M))) &\stackrel{\text{Def}}{=}
  KO_0(\Cliff_n\hat\tensor
  C^*_{red}(\gamma))\\
 &=KO_0(\Cliff_{n+1}\hat\tensor
  C^*(\Z/2\times\Z^3)) = KO_{n+1}(C^*(\Z/2\times\Z^3)) .
\end{split}\end{equation*}
For a discrete group $G$ we have
(cf. \cite[p.~14 and 1.5.4]{Schroeder(1993)})
\begin{equation*}
KO_n(C_{red}^*(G\times\Z))\iso KO_n(C_{red}^*(G))\oplus
KO_{n-1}(C_{red}^*(G)).
\end{equation*}
Therefore the lemma now follows from the table below. We compute it
using the fact that $\integers/2$
has exactly two irreducible real representation, both of real
type. Hence $KO_*(\reals\integers/2)=KO_*(\reals)\oplus KO_*(\reals)$.\\

\begin{tabular}{l|l|l|l|l|l|l|l|l}
$n \pmod 8$ & $1 $ & $ 2 $ & $ 3 $ & $ 4 $ & $ 5 $ & $ 6 $ & $ 7 $ & $
0$\\ \hline
$KO_n(C^*\Z/2) $ & $ \Z/2\times \Z/2 $ & $ \Z/2\times \Z/2 $ & $ 0 $ & $\Z^2 $ & $ 0 $ & $ 0 $ & $
0 $ & $ \Z^2$
\end{tabular}

\end{proof}

\begin{corollary}
  The manifold $M$ defined in Section \ref{sec_ex} has twisted index
  \begin{equation*}
    \alpha([M]) = 0 \in KO_{5}(C^*_{red}(\gamma(M))).
  \end{equation*}
\end{corollary}

\begin{proof}
  We have $\Omega_{*}(\gamma(M)) \cong \Omega_{*}^{pin'}(B\pi')$. Since the 
  coefficients of $Pin'$-bordism theory are completely $2$-torsion 
(cf. \cite{Giambalvo(1973)}) it follows from the Atiyah-Hirzebruch 
spectral sequence that $\Omega_{*}(\gamma(M))$ is completely 2-torsion. 
Hence $\alpha: \Omega_{5}(\gamma(M)) \to  KO_5(C^*_{red}(\gamma)) =\Z^6$ 
is trivial.
\end{proof}

%%% Local Variables: 
%%% mode: latex
%%% TeX-master: "GLR.Michael"
%%% TeX-master: "GLR.Michael"
%%% TeX-master: "januar_version"
%%% TeX-master: "januar_version"
%%% TeX-master: "februar_version"
%%% End: 

\end{document}